\documentclass[12pt]{article}
\usepackage[utf8]{inputenc}
\pdfoutput=1 
\usepackage{amsfonts}
\usepackage{amsmath}
\usepackage{amssymb}
\usepackage{amsthm}

\usepackage[letterpaper, margin=3cm]{geometry}

\usepackage{graphicx} 
\usepackage[hang]{caption}

\usepackage[numbers, square,sort&compress,comma]{natbib}

\usepackage[bookmarks=false, colorlinks=true, linkcolor=blue, urlcolor=blue, citecolor=blue]{hyperref}
\usepackage[framestyle=fbox,heightadjust=all,framearound=all]{floatrow}

\usepackage{pifont}
\usepackage{setspace}
\usepackage[titletoc, title]{appendix}

\begin{document}
\title{A second-order, perfectly matched layer formulation to model 3D transient wave propagation in anisotropic elastic media}

\author{Hisham Assi\footnote{\textit{Email address:} hisham.assi@mail.utoronto.ca}\quad and Richard S. C. Cobbold\\ \small {Institute of Biomaterials and Biomedical Engineering, University of Toronto,}\\ \small{164 College Street, Toronto, M5S 3G9, Canada}}
\date{}
\maketitle


\begin{abstract}
Numerical simulation of wave propagation in an infinite medium is made possible by  surrounding a finite region by a perfectly matched layer (PML). Using this approach a generalized three-dimensional (3D) formulation is proposed for time-domain modeling of elastic wave propagation in an unbounded lossless anisotropic medium. The formulation is based on a second-order approach that has the advantages of, physical relationship to the underlying equations, and amenability to be implemented in common numerical schemes. Specifically, our formulation uses three second-order equations of the displacement field  and nine auxiliary equations, along with the three time histories of the displacement field. The properties of the PML, which are controlled by a complex two-parameter stretch function, are such that it acts as  near perfect absorber.  Using finite element method (FEM) 3D numerical results are presented for a highly anisotropic medium.  An extension of the formulation  to the particular case of a Kelvin-Vogit viscoelastic medium is also presented. 
\vspace{5mm}\\ 
\textit{Keyword}: Perfectly matched layers; Elastic waves; Second order time-domain; Anisotropic media; viscoelastic media 

\end{abstract}

 \tableofcontents
\doublespacing
\section{Introduction}
\label{sec:introduction}
Perfectly matched layers \citep{Berenger:1994ua} are a well-developed method for simulating wave propagation in unbounded media enabling the use of a reduced computational domain without having to worry about spurious boundary reflections. B\'{e}renger  showed that by adding specific conductivity parameters to Maxwell's equations, perfect matching and decaying of the propagating waves in the PML could be achieved \citep{Berenger:1994ua}. An alternative method is to assume that the material contained within the PML is a uniaxial anisotropic media \citep{Sacks:1995gs,Roden:1997fo,Gedney:1996ub}, generally referred to as the uniaxial PML approach. A third method, with greater generality and flexibility, is the complex coordinate stretching approach \citep{Chew:1994dn}. In fact, the conductivity parameter introduced by B\'{e}renger \citep{Berenger:1994ua} can be thought of as a damping parameter in a stretch function that transforms the spatial coordinate in the layer to the complex plane.
 
Subsequent to these electromagnetic wave applications, many PML formulations have been introduced for elastic wave propagation \citep{Collino:2001vt,Hastings:1996um,Drossaert:2007fi,Chew:1996wk,Appelo:2006wz,MezaFajardo:2008dx,Kucukcoban:2011tn}.  Amongst these, the split-field formulations that are typically described by systems of first order equations with double (for 2D) or triple (for 3D) the number of stress-velocity physical equations (nine equations in all for the 3D case) . Second-order formulations uses one physical field variable (usually the displacement) along with extra auxiliary variables that are typically needed to obtain the time-domain equations from the frequency-domain equations. Mathematically, it has been proven that certain second-order  PMLs are strongly well-posed, while the first-order type is only weakly well-posed \citep{Deinega:2011er,Abarbanel:2002wm,Duru:2012va}.

There are other advantages for choosing second-order formulations. The second-order displacement elastic wave equation emerges directly from Newton's second law \citep {Duru:2012va}, unlike the fist-order stress-velocity elastic wave equation which introduces a new non-physical wave mode with zero velocity \citep{Duru:2012va, Appelo:2006wz}. Moreover, the second-order PML formulations are more readily implemented in common numerical schemes that are based on second-order displacement wave equations \citep{Komatitsch:2007bz, Li:2010tu}. However, deriving  second-order PML formulations is less trivial than that for first-order ones, especially in the time-domain where many auxiliary variables are needed. The problem becomes more complex for the 3D modeling which would partially explain the dearth of second order formulations in 3D.

There have been a number of papers that describe the formulation of time-domain wave propagation in 3D fluid media using PMLs (see of example \citep{Hu:2007jb,kaltenbacher2013modified}) but there  are relatively few that address the same problem for anisotropic, inhomogeneous  elastic media, especially those attempting a second-order formulation. In  previous works \citep{Assi:2013tp, Assi:2015tv}, the authors introduced a compact second-order time-domain PML formulations for the elastic wave equation in 2D which has only four auxiliary variable. Recently Lee and Shin \citep{Koria:2015th} introduced an unsplit PML formulation for isotropic media or media with vertical axis of symmetry (VTI).  Their 2D derivation was based on second-order elastic wave equations, and the final formulation followed closely the one given in Assi and Cobbold \citep{Assi:2013tp, Assi:2015tv}. It should be noted that the final form of the PML formulation, and not the way it wavs derived, that governs its robustness and other characteristics. Additionally, Lee and Shin  extended their formulation to 3D VTI media for which they presented numerical results \citep{Koria:2015th}. 

The purpose of this paper is to derive a time-domain second-order formulation to model elastic wave propagation in an unbounded three-dimensional general anisotropic inhomogeneous  solid.  As will be seen the  formulation results in a system of equations that are applicable throughout the computational domain. In the physical domain, the complex stretch function is simply set to unity.  To demonstrate the application of our formulation, propagation from a spherical transient source  embedded in a highly anisotropic medium (the mineral Olivine) is illustrated. Extension of the formulation to include a viscoelastic medium that can be represented by a Kelvin-Vogit model \citep{Meyers:2008,Banks:2011gj}, is presented in \autoref{app:viscoelasticity}.

\section{Background and materials}
\label{sec:background}
\subsection{Elastic waves in solids}
\label{subsec:theoElastic}
Wave propagation in linear elastic solids can be described using Newton's second law, along with Hook's law and a linear approximation for the strain. These lead to the following second-order formulation of the elastic wave equation: 

\begin{equation}
\rho\frac{\partial^{2}u_{i}}{\partial t^{2}}-\sum\limits _{j=1}^{3}\frac{\partial}{\partial x_{j}}\left(\sum\limits _{k,l=1}^{3}C_{ijkl}\frac{\partial u_k}{\partial x_l}\right)=0,
\label{eq:waveElastic}
\end{equation}
where $t\in\mathbb{R}^{+} $ is time, $\mathbf{x}\in\mathbb{R}^{3}$ is the space variable, $u_i(\mathbf{x},t)$ are the components of particle displacement vector. Moreover, $\rho(\mathbf{x})$ is the solid mass density and $C_{ijkl}(\mathbf{x})$ are the components of the fourth order elasticity tensor with the following symmetry properties: $C_{ijkl}=C_{ijlk}=C_{jikl}$, and $C_{ijkl}=C_{klij}$. The source of energy that excites the elastic medium can be added as a load vector, $F_i(\mathbf{x},t)$, to the right-hand side (RHS) of \eqref{eq:waveElastic}.

In general, the elasticity tensor, $C_{ijkl}$, has 81 components, but due to the above symmetries, the maximum number of independent parameters is 21. For the special case of isotropic solids, the elasticity tensor can be described by two independent parameters such as the  Lam\'{e} coefficients,  $\lambda(\mathbf{x})$ and $\mu(\mathbf{x})$. In terms of these two coefficients, the elasticity tensor can be written as:
\begin{equation}
C_{ijkl}=\lambda\delta_{ij}\delta_{kl}+\mu(\delta_{ik}\delta_{jl}+\delta_{il}\delta_{jk}),
\label{eq:HookIsotropic}
\end{equation}
where $\delta_{ij}$is the Kronecker delta function.

For the purpose of plane wave and Fourier analyses, the harmonic wave solutions of the following form:
\begin{equation}
\mathbf{u}=\mathbf{A}\exp{\left[i(\mathbf{k\cdot x}-\omega t)\right]},
\label{eq:harmonicWave}
\end{equation}
will be considered for the elastic wave equation as given by \eqref{eq:waveElastic}. In this equation, $\mathrm{\mathbf{A}} \in \mathbb{C}^{3}$ 
is the constant amplitude polarization vector, $\mathrm{\mathbf{k}}\in\mathbb{R}^{3}$ is the wavevector, $\omega\in\mathbb{C}$ is the angular frequency, and $i^2=-1$. 

\subsection{Complex coordinates stretching}
\label{subsec:complex}
\sloppy{To obtain a PML formulation for a given wave equation, the complex coordinate stretching \citep{Chew:1994dn} can expressed as a coordinate transform: ${\mathbf{x}\mapsto\tilde{\mathbf{x}}:\mathbb{R}^3\to\mathbb{C}^3}$.}  Since  $\tilde{\mathbf{x}}=\mathbf{x}$ in  the physical domain and the PML region is assumed to be homogeneous, then $\tilde{\mathbf{x}}$ appears only in the form of spatial partial derivatives in the PDEs. Given a field variable $u$, then using the chain rule: $\frac{\partial u} {\partial x_j}= \sum_{k=1}^{3} \frac{\partial u}{\partial \tilde{x}_k} \;\frac{\partial \tilde{x}_k}{\partial x_j}$, which reduces to $\frac{\partial u} {\partial x_j}= \frac{\partial u}{\partial \tilde{x}_j} \;\frac{\partial \tilde{x}_j}{\partial x_j}$ since $\tilde{x}_j$ depends only on $x_j$.  As a result, defining the complex stretch function by $s_{j}\left(x_{j}\right)=\frac{\partial\tilde{x}_{j}}{\partial x_{j}}$ suffices to perform transformation:
\begin{equation}
\frac{\partial}{{\partial\tilde{x}}_{j}}=\frac{1}{s_{j}\left(x_{j}\right)}\frac{\partial}{{\partial x}_{j}}.
\label{eq:derivativeTransform}
\end{equation}

The two-parameter complex stretch function introduced by Fang and Wu \citep{Fang:1996tc} in their generalized PML (GPML)  is adopted in this paper. This function is given by
\begin{equation}
s_{j}\left(x_{j}\right)=\alpha_{j}\left(x_{j}\right)\left[1+i\,\frac{\beta_{j}\left(x_{j}\right)}{\omega}\right],
\label{eq:CSF}
\end{equation}
where the $\beta_{j}\geq0$ is the damping parameter  responsible for damping the propagating wave inside the PML. In this equation, the scaling parameter, $\alpha_{j}>0$, is responsible for either stretching ($\alpha_{j}>1$) or compressing ($0<\alpha_j<1$) the coordinate. It should be noted that in the physical domain, where $\tilde{x}_{j}\left(x_{j}\right)=x_{j}$, $\beta_{j}=0$ and $\alpha_{j}=1$.

Appropriate choices are now needed for the stretch function parameters $\alpha_{j}(x_{j})$ and $\beta_{j}(x_{j})$. Despite the absence of a rigorous methodology for their choice \citep{Chew:1996wk,Kucukcoban:2011tn}, polynomial functions are often used as indicated below for the damping parameter:
\begin{equation}
%
\beta_{j}\left(x_{j}\right)  =\begin{cases}
0 & \text{if}\left|x_{j}\right|<x_{0}\\
\beta_{0_{j}}\left(\dfrac{\left|x_{j}\right|-x_{0}}{d}\right)^{n} & \text{if }x_{0}\leq\left|x_{j}\right|\leq x_{0}+d,
\end{cases}
\label{eq:beta}
\end{equation}
$\!\!$where \textit{d} is the thickness of the PML, ${2x}_{0}$ is the dimension of the square physical domain centered at the origin, $n$ is the polynomial order,  $\beta_{0_{j}}$ is a constant that represent the maximum values of $\beta_{j}$. The value of this parameter needs to be specified. It is helpful to express the value of $\beta_{0_{j}}$ in terms of a desired amplitude reflection coefficient ($R_{j})$ caused by the reflection from the outer boundary of the PML. It can be shown that for normal incidence and assuming $\alpha_{j}=1$,
\begin{equation}h
\beta_{0_{j}}=\frac{c_{\mathrm{max}}\left(n+1\right)\ln\left(1/R_{j}\right)}{2d}.
\label{eq:reflection}
\end{equation}
Quadratic polynomial, corresponds to $n=2$, will be used in this work unless mentioned otherwise.
Without loss of generality, the scaling parameter,  $\alpha_{j}$, is set to unity in this work. This parameter can be readily introduced back to any PML formulation that is derived using the complex stretch function in \eqref{eq:CSF}, by replacing each $\partial/\partial x_j$ by $\partial/ \alpha_j\partial x_j$ in the PDEs.
\section{Formulation of PML for elastic wave propagation}
\label{sec:formulation} 
With the help of the above background, our time-domain PML formulation can be introduced for the  wave propagation in unbounded media. All parameters, namely, $\rho$, $C_{ijkl}$, $\beta_{j}$, and $s_{j}$, are assumed to be space dependent throughout the derivation leading to a variable-coefficient PML formulation. Since the stretch function also depends on the frequency, the derivation starts in the frequency domain. 

First, we take Fourier transforms of the elastic wave equation \eqref{eq:waveElastic}, and then transform the spatial coordinates using complex stretching, ${\mathbf{x}\mapsto\tilde{\mathbf{x}}}$, as introduced in \autoref{subsec:complex}. These steps lead to:
\begin{equation}
\left(-i\omega\right)^{2}\hat{u}_i\,\rho=\sum\limits _{j=1}^{3}\frac{\partial}{\partial \tilde{x}_{j}}\left(\sum\limits _{k,l=1}^{3}C_{ijkl}\frac{\partial \hat{u}_k}{\partial \tilde{x}_l}\right),
\label{eq:derevationElastic1}
\end{equation}
where $\hat{\Box}$ denotes the Fourier transform in time. The  need to solve this differential equation along contours in the complex plane can  be avoided by inverse transforming the complex-stretched coordinates back to the original spatial coordinates using \eqref{eq:derivativeTransform}. This is followed by multiplying the equation by  $s_1\, s_2 \, s_3$, leading to:
\begin{equation}
s_1s_2s_3\,\left(-i\omega\right)^{2}\hat{u}_i\,\rho=\sum\limits _{j=1}^{3}\frac{s_1s_2s_3}{s_j}\frac{\partial}{\partial x_{j}}\left(\sum\limits _{k,l=1}^{3}C_{ijkl}\frac{1}{s_l}\frac{\partial \hat{u}_k}{\partial x_l}\right).
\label{eq:derevationElastic2}
\end{equation}
Expanding  $s_1s_2s_3$ according to \eqref{eq:CSF} while assuming $\alpha_j=1$, the left-hand side (LHS) of the above equation becomes:
\begin{equation}
\rho\,\left[ \left(-i\omega\right)^{2}+(-i\omega) \left(\beta_1+\beta_2+\beta_3\right)+ \left(\beta_1\,\beta_2+\beta_2\,\beta_3+\beta_3\,\beta_1 \right)+\frac{\beta_1\beta_2\beta_3}{(-i\omega)} \right] \,\hat{u}_i.
\label{eq:derevationElastic31}
\end{equation}
Here, it is helpful to introduce the variable
\begin{equation}
U_i(\mathbf{x},t)=\int_{0}^{t}u_i(\mathbf{x},\tau)\, d\tau,
\label{eq:derevationHistory}
\end{equation}
 whose Fourier transform is given by $\hat{U}_j\mathbf{\mathrm{(}x},\omega)=\hat{u}_i(\mathbf{x},\omega)/(-i\omega)+\pi\,\hat{u}_i(\mathbf{x},0)\delta(\omega)$. However, the second term vanishes since the stretch function \eqref{eq:CSF} is not defined for the static case of $\omega=0$ \citep{Kucukcoban:2011tn}. Consequently, substituting $\hat{u}_{i}=-i\omega\,\hat{U}_{i}$ in the last term of \eqref{eq:derevationElastic31} and taking  inverse Fourier transform $(-i\omega\Rightarrow\thinspace\partial/\partial t)$ of this, results in
\begin{equation}
\rho\,\left[ \frac{\partial^{2}u_i}{\partial t^2}+\left(\beta_1+\beta_2+\beta_3\right)\frac{\partial u_i}{\partial t}+ \left(\beta_1\,\beta_2+\beta_2\,\beta_3+\beta_3\,\beta_1 \right)\,u_i+\beta_1\beta_2\beta_3\,U_i\right].
\label{eq:derevationElastic32}
\end{equation} 

It should be noted that $s_1(x_1)s_2(x_2)s_3(x_3)/s_j(x_j)=\prod_{i\neq j} \,s_i(x_i)$, and hence does not depend on $x_j$, enabling  the RHS of \eqref{eq:derevationElastic2} to be written as:
\begin{equation}
\sum\limits _{j=1}^{3}\frac{\partial}{\partial x_{j}}\left(\sum\limits _{k,l=1}^{3}\frac{s_1s_2s_3}{s_j s_l}\,C_{ijkl}\,\frac{\partial \hat{u}_k}{\partial x_l}\right).
\label{eq:derevationElastic4}
\end{equation}
After some manipulations, it can be shown that
\begin{equation}
\frac{s_1s_2s_3}{s_j s_l}=1+\frac{\beta_1+\beta_2+\beta_3-\beta_j-\beta_l+\frac{\beta_1\beta_2\beta_3}{-i\omega\,\beta_l}}{-i\omega+\beta_j}
\label{eq:derevationElastic41}
\end{equation}
At this point, we introduce the auxiliary variables, $w_{ij}(\mathbf{x},t)$ such that their Fourier transform
\begin{equation}
\hat{w}_{ij}(\mathbf{x},\omega)=\sum\limits _{k,l=1}^{3} \frac{\beta_1+\beta_2+\beta_3-\beta_j-\beta_l+\frac{\beta_1\beta_2\beta_3}{-i\omega\,\beta_l}}{-i\omega+\beta_j}\;C_{ijkl}\frac{\partial \hat{u}_k}{\partial x_l}.
\label{eq:derevationElastic42}
\end{equation}
Multiplying the above equations by $-i\omega+\beta_j$ and taking its inverse Fourier transform results in the following time-domain auxiliary equations:
\begin{equation}
\frac{\partial w_{ij}}{\partial t}+\beta_{j}\, w_{ij}  =\sum_{k,l=1}^3\left(\left(\beta_1+\beta_2+\beta_3-\beta_j-\beta_l\right)\,C_{ijkl}\,\frac{\partial u_k}{\partial x_l}+ \frac{\beta_1\beta_2\beta_3}{\beta_l}C_{ijkl}\,\frac{\partial U_k}{\partial x_l}\right),
\label{eq:derevationAuxiliary}
\end{equation}
and the RHS of \eqref{eq:derevationElastic2} becomes
\begin{equation}
\sum\limits _{j=1}^3\frac{\partial}{\partial x_j}\left(\sum\limits _{k,l=1}^3 C_{ijkl}\frac{\partial u_k}{\partial x_l}+w_{ij}\right).
\label{eq:derevationElastic43}
\end{equation}
 
This concludes our derivation, so that the final second-order time-domain PML formulation for elastic wave propagation in three-dimensional anisotropic solid can be written as
\begin{subequations}
\begin{align}
\rho\left(\frac{\partial^{2}u_i}{\partial t^2}+a\,\frac{\partial u_i}{\partial t}+b\, u_i +c\, U_i\right) & =\sum\limits _{j=1}^3\frac{\partial}{\partial x_j}\left(\sum\limits _{k,l=1}^3 C_{ijkl}\frac{\partial u_k}{\partial x_l}+w_{ij}\right)\\[12pt]
\frac{\partial w_{ij}}{\partial t}+\beta_{j}\, w_{ij} & =\sum_{k,l=1}^3\left(\tilde{C}_{ijkl}\,\frac{\partial u_k}{\partial x_l}+\breve{C}_{ijkl}\,\frac{\partial U_k}{\partial x_l}\right)\\[12pt]
\frac{\partial U_i}{\partial t}&=u_i ,
\end{align}
\label{eq:pmlElastic}
\end{subequations}
where $a(\mathbf{x})=\beta_1+\beta_2+\beta_3$, $b(\mathbf{x})=\beta_1\,\beta_2+\beta_2\,\beta_3+\beta_3\,\beta_1$, $c(\mathbf{x})=\beta_1\,\beta_2\,\beta_3$, $\tilde{C}_{ijkl}(\mathbf{x})=(a-\beta_j-\beta_l)\,C_{ijkl}$, and $\breve{C}_{ijkl}(\mathbf{x})=(c \,/\, \beta_l)\, C_{ijkl}$.

\section{Numerical Methods and Results}
\label{sec:numerical} 
For our studies, the source of excitation was a 1~mm radius sphere, centered at the origin and embedded in an infinite 3D medium. To model the infinite medium we assumed a cubic physical domain of 2.0 cm$^{3}$ that is centered at the origin and surrounded by a 2.0~mm PML. The boundary of the sphere was assumed to vibrate  with a displacement, whose normalized time-dependence is given by the first derivative of a Gaussian, i.e., by
\begin{equation}
u_{0}\left(t\right)=-\sqrt{2e}\thinspace\pi f_{0}\left(t-t_{0}\right)\thinspace e^{-\pi^{2}f_{0}^{2}\left(t-t_{0}\right)},
\label{eq:source}
\end{equation}
where $f_{0}$ is the dominant frequency and $t_{0}$ is a source delay time. All numerical experiments used $f_{0}=1$ MHz and $t_{0}=1~\mu$s.

\begin{figure}[H]
\centering
\mbox{\includegraphics[width=1.0\textwidth]{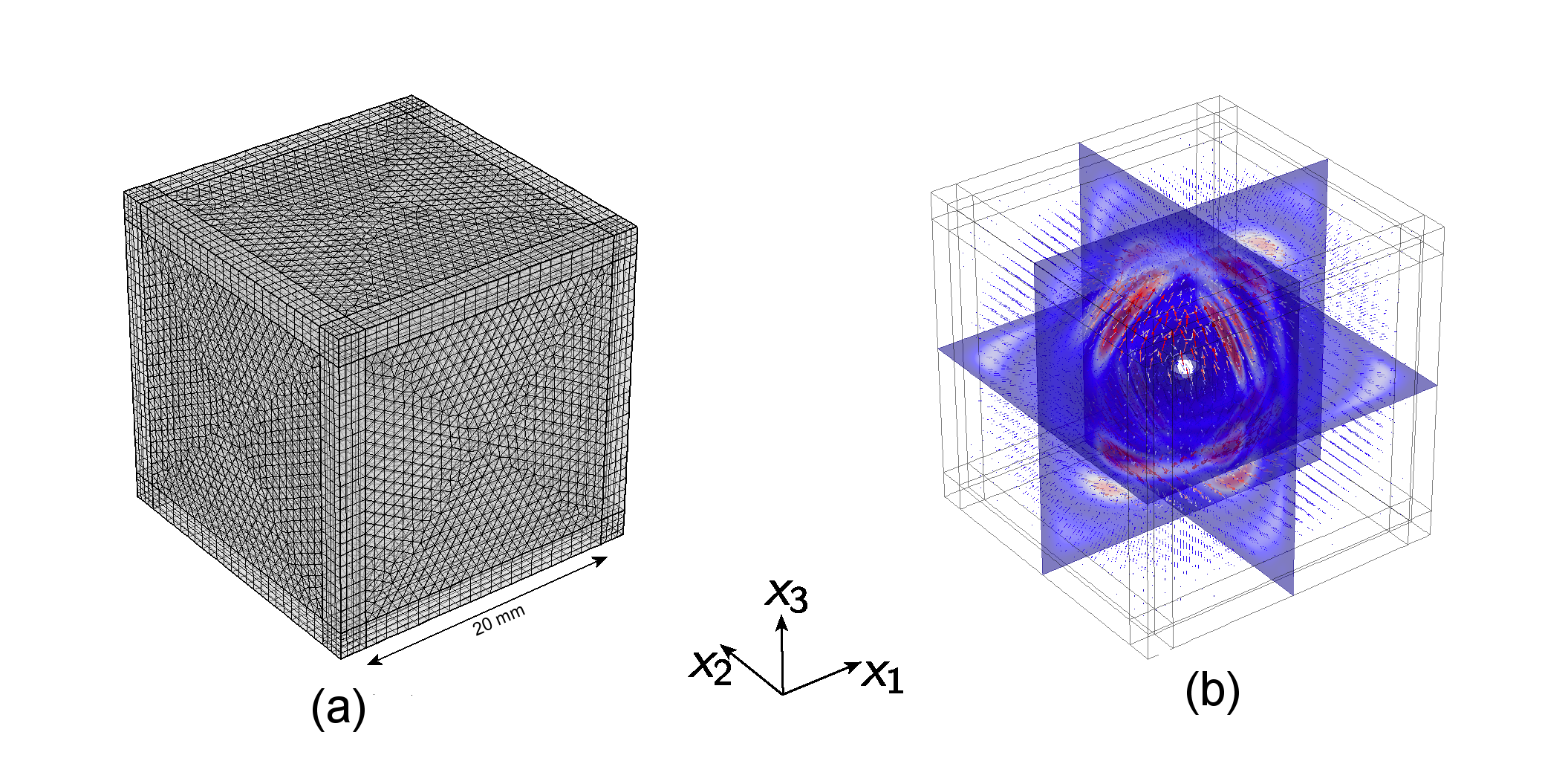}}
\protect\caption{(a) Showing the 3D grid used in numerical calculations. (b) Snapshot at $t=2.5~\mu$s.illustrating the 3D response for the mineral Olivine.}
\label{fig:ThreeDGeometry}
\end{figure}

The simulations were performed using COMSOL with the second-order Lagrange finite elements employing a cubic mesh for the PML region and the default  tetrahedral  shape in the physical domain (see \autoref{fig:ThreeDGeometry}a). For discussing the mesh dimensions and time discretization, it is  helpful to define the minimum and maximum characteristic wave speeds associated with the medium by $c_{\text{min}}$ and $c_{\text{max}}$.  The mesh size is governed by the shortest wavelength of significance for the propagating pulse, i.e., by $c_{\text{min}}/f_c$.
 Specifically, the mesh size was assumed to be given by
\begin{equation}
h_{0}=\frac{c_{\text{min}}}{Nf_{0}},\label{eq:meshSize}
\end{equation}
which for  the  second-order accurate  finite elements  used, corresponds to $2N$ degrees of freedom per wavelength. As illustrated in \autoref{fig:ThreeDGeometry}a, the mesh employed uses a PML whose thickness consists of just four elements.  As will be seen this is sufficient to ensure virtually complete absorption of the various incident waves. For the time discretization, a second-order generalized alpha method, as defined by Chung and Hulbert \citep{Chung:1993io}, was used with $\rho_{\infty}$= 0.75. The step duration was $0.9h_{0}/c_{\text{max}}$, which is just less than the time needed for the fastest wave to travel through the smallest mesh dimension.

To test the validity of our formulation and the accuracy with which our finite element simulations describe the propagating pulse, the exact solution for a monochromatic compressional wave caused by a 1-mm radius sphere whose surface vibrated normal to the surface \citep{Beltzer:1988ut} was used. The sphere was assumed to be embedded in an unbounded isotropic solid (glass). By multiplying this monochromatic solution with the Fourier transform of equation \eqref{eq:source} and then taking the inverse Fourier transform, the time-domain analytic solution was obtained. Good agreement with the numerically calculated response  provided convincing evidence for the validity of our 3D formulation, though the results are not presented here.

To illustrate  the application of our general formulation, we chose to present the results for a highly anisotropic medium. In particular, we chose  to examine the 3D response when the above source, in \eqref{eq:source}, is contained within in a single crystal of olivine (Mn$_{2}$SiO$_{4}$). Olivine is a mineral with an orthorhombic structure and nine independent elasticity components whose measured parameters at 25$^{\circ}$C are given by \citep{sumino1979elastic,bass1984elasticity}: $c_{11}=2.58$, $c_{22}=1.66$, $c_{33}=2.07$, $c_{44}=0.45$, $c_{55}=0.56$, $c_{66}=0.58$, $c_{12}=0.87$, $c_{13}=0.95$, $c_{23}=0.92$ Mbar. Unlike the isotropic case,  fast and slow waves propagate in anisotropic solids even if the excitation is normal to the sphere's surface. Nevertheless, in order to observe a clearer presence of these different waves, we decided to to excite the medium by vibrating the sphere's surface at $45^\circ$ to the normal in the polar direction. Namely, the Dirichlet boundary condition at the surface of the sphere is set to:
\begin{equation}
\frac{1}{\sqrt{2}}(\hat{\mathbf{n}}+\hat{\mathbf{t}}_{\phi})\,u_{0}\left(t\right),
\label{eq:Dirichlet}
\end{equation}
where $u_0$ is defined in \eqref{eq:source}, $\hat{\mathbf{n}}$ is the normal unit vector, and $\hat{\mathbf{t}}_{\phi}$ is the tangential unit vector in the $\phi$ direction, and $\phi$ is the polar angle that varies from 0 to $\pi$ away from the $x_3$-axis. At $\phi=90^\circ$ for example,  $\hat{\mathbf{t}}_{\phi}$ is in the negative $x_3$-direction, hence, the quasi-longitudinal wave is expected to be dominant on the $x_1$-$x_2$~plane.

\begin{figure}[H] 
\centering
\mbox{ \includegraphics[width=1.0\textwidth]{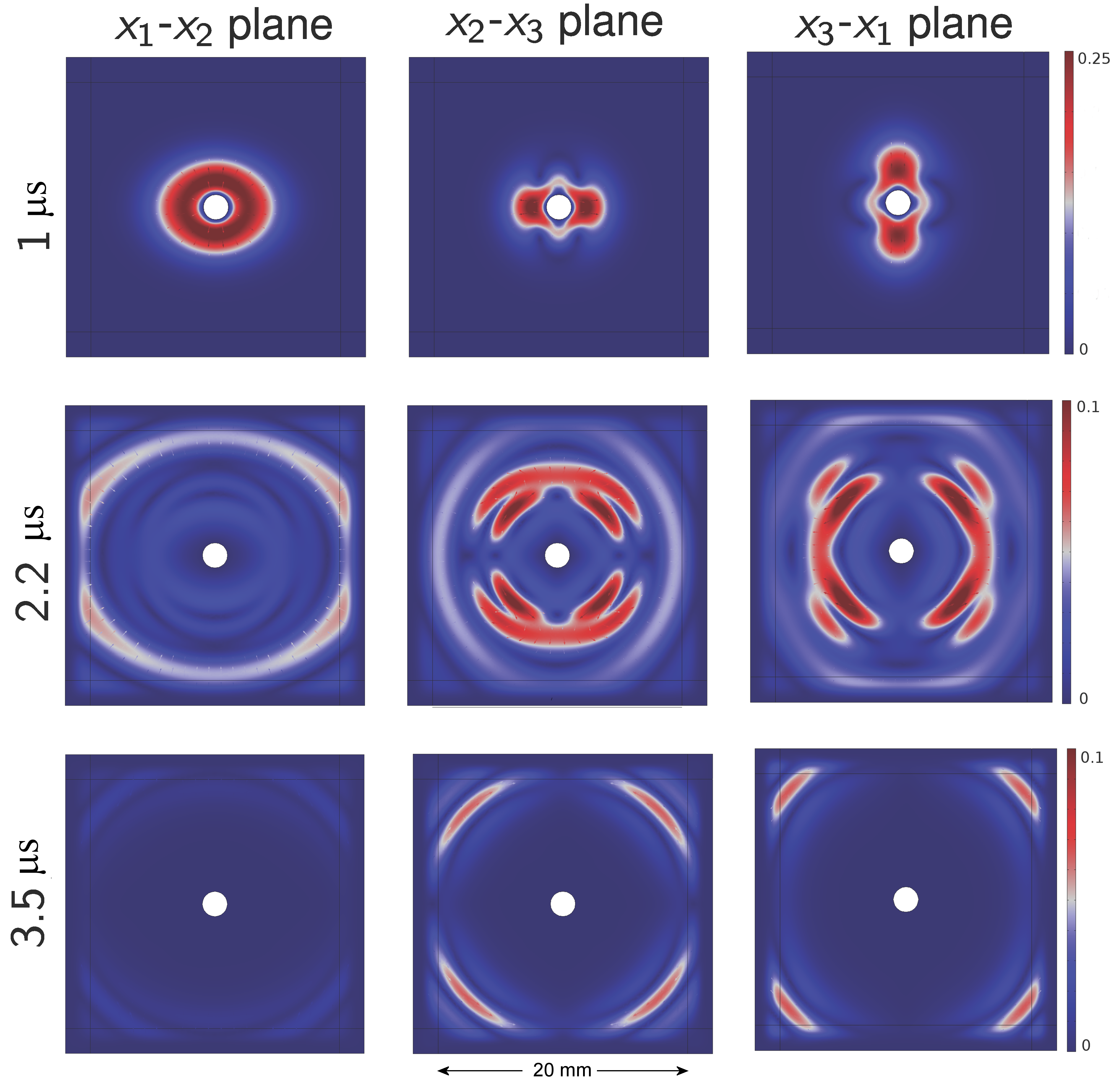}}
\caption{Sets of snapshots of the waves produced by the 1-mm radius source excited by the waveform given by equation \eqref{eq:source},  at three different times for three different planes.}
\label{fig:Planes2D}
\end{figure}

The results of the simulations are presented as density and vectors plots that represent the magnitudes and the directions of the normalized displacement field. While \autoref{fig:ThreeDGeometry}b provides a snapshot of the propagating waves in a 3D format,  such an image is difficult to interpret. The three sets of snapshots for three different planes, as shown in the nine panels of  \autoref{fig:Planes2D}, provides much more detailed information. These snapshots show 2D plots of the field on the thee principal planes at 1~$\mu$s, 2.2~$\mu$s, and 3.5~$\mu$s. The first column shows the displacement field on the $x_1$-$x_2$~plane, wherein, as expected, the fast wave is dominant. Meanwhile,  on the other two planes, the thee waves; the quasi-longitudinal (fast) and the two quasi-shear (slow), are clearly present as shown in the 2.2~$\mu$s snapshots. At this time, the fast wave is being effectively absorbed by the PML, while the slow waves are being absorbed  in the  3.5~$\mu$s snapshot.   

\begin{figure}[H]
\centering
\mbox{\includegraphics[width=0.7\textwidth]{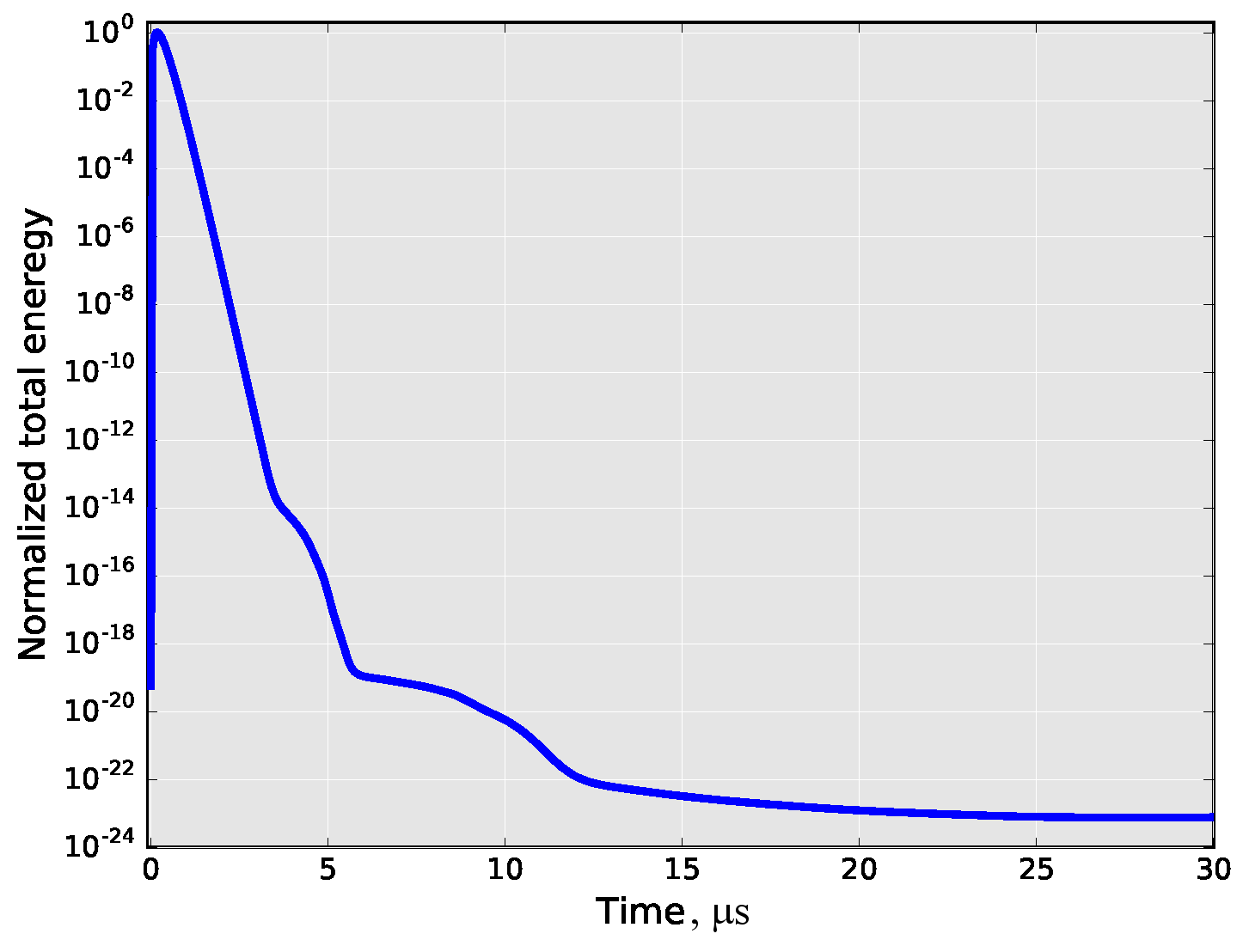}}
\protect\caption{The decay of energy as calculated from \eqref{eq:TotalEnergy}  in the physical domain for the mineral Olivine showing that the PML acts as a near ideal absorber.}
\label{fig:EnergyDecay}
\end{figure}

The effectiveness of the PML to absorb all the incident energy can be obtained by looking at the manner in which the energy in the physical domain evolves over time. Since the total energy in the physical domain is the sum of the kinetic and potential energy, it can be calculated from
\begin{equation}
E (t)=\frac{1}{2}\int\limits _{\Omega}\left[\rho\sum\limits _{j=1}^{3}\left(\dfrac{\partial u_{j}}{\partial t}\right)^{2}+\sum\limits _{i,j,k,l=1}^{3}C_{ijkl}\dfrac{\partial u_{i}}{\partial x_{j}}\dfrac{\partial u_{k}}{\partial x_{l}}\right]\mathrm{d}\Omega.\label{eq:TotalEnergy}
\end{equation}
\autoref{fig:EnergyDecay} shows that the total energy decays to a negligible level in less than 5$~\mu$s.

\section{Conclusion}
\label{sec:conclusion} 
Using PML approach we have addressed the problem of wave propagation in an unbounded, linear anisotropic solid in three dimensions. A time-domain second order PDE has been derived using complex coordinate stretching. The advantages of our time-domain formulation is the fact that it covers the more general inhomogeneous anisotropic case using a small number of equations. Specifically, three second-order equations of the displacement field  and nine auxiliary equations, along with the time histories of the displacement field. This simplifies the problem and reduces the computational resources needed. Moreover,  use can be made of a wider variety of second-order numerical schemes.

\section*{Acknowledgments}
We wish to thank the Natural Sciences and Engineering Research Council (NSERC) for their support [grant number 3247-2012]. We also wish to thank Pooya Bidari from Ryerson University for drawing to our attention the importance of developing a viscoelastic 3D model.

\newpage
\bibliographystyle{IEEEtranN}
\bibliography{ThreeD_2}

\newpage
\begin{appendices}
\section{PML for viscoelastic media}
\label{app:viscoelasticity}
It should be noted that the PML formulation presented in the this paper models  wave propagation in a loss-less media, for which, Hook's law as used in \eqref{eq:waveElastic}, is given by
\begin{equation}
\sigma_{ij}=\sum\limits _{k,l=1}^{3}C_{ijkl}\frac{\partial u_k}{\partial x_l},
\label{eq:Hooks}
\end{equation}
where $\sigma_{ij}$ are the components of the symmetric stress tensor, and the components of the elasticity tensor, $C_{ijkl}$, are assumed to be real-valued. There are several different models that are used to account for viscous losses to the elastic wave equation \citep{Banks:2011gj}. One that is commonly used especially for modeling  wave propagation in tissue \citep{Bercoff:2004es},   is the Kelvin--Vogit model, for which Hook's law takes the form
\begin{equation}
\sigma_{ij}=\sum\limits _{k,l=1}^{3}\left(C_{ijkl}+ \eta_{ijkl}\,\frac{\partial}{\partial t} \right) \frac{\partial u_k}{\partial x_l},
\label{eq:Hooks}
\end{equation}
where $\eta_{ijkl}$ is the viscosity tensor. For such a medium the wave equation is:
\begin{equation}
\rho\frac{\partial^{2}u_{i}}{\partial t^{2}}=\sum\limits _{j=1}^{3}\frac{\partial}{\partial x_{j}}\left(\sum\limits _{k,l=1}^{3}C_{ijkl}\frac{\partial u_k}{\partial x_l}+ \eta_{ijkl}\frac{\partial^2 u_k}{\partial t\,\partial x_l}\right).
\label{eq:waveVescoElastic}
\end{equation}

Following the same steps used in deriving the elastic PML formulation in \autoref{sec:formulation}, a PML formulation for the above viscoelastic wave equation can be obtained. The only difference, for this viscoelastic case, is that equation~\eqref{eq:derevationElastic4} in the derivation becomes:
\begin{equation}
\sum\limits _{j=1}^{3}\frac{\partial}{\partial x_{j}}\left(\sum\limits _{k,l=1}^{3}\frac{s_1s_2s_3}{s_j s_l}\,\left(C_{ijkl}-i\omega \,\eta_{ijkl}\right) \,\frac{\partial \hat{u}_k}{\partial x_l}\right).
\label{eq:derevationVisco}
\end{equation}
This leads to the following PML formulation for the viscoelastic  wave equation:
\begin{subequations}
\small
\begin{flalign}
\rho_S\left(\frac{\partial^{2}u_i}{\partial t^2}+a\,\frac{\partial u_i}{\partial t}+b\, u_i +c\, U_i\right) =\sum\limits _{j=1}^3\frac{\partial}{\partial x_j}\left(\sum\limits _{k,l=1}^3 C_{ijkl}\frac{\partial u_k}{\partial x_l}+ \eta_{ijkl}\frac{\partial^2 u_k}{\partial t\,\partial x_l}+w_{ij}\right)\\[12pt]
\frac{\partial w_{ij}}{\partial t}+\beta_{j}\, w_{ij}=\sum_{k,l=1}^3\left(\tilde{\eta}_{ijkl}\frac{\partial^2 u_k}{\partial t\,\partial x_l}+ (\tilde{C}_{ijkl}+\breve{\eta}_{ijkl}) \,\frac{\partial u_k}{\partial x_l}+\breve{C}_{ijkl}\,\frac{\partial U_k}{\partial x_l}\right)\\[12pt]
\frac{\partial U_i}{\partial t}=u_i ,
\end{flalign}
\label{eq:pmlVElastic}
\end{subequations}
where $\tilde{\eta}_{ijkl}(\mathbf{x})=(a-\beta_j-\beta_l)\,\eta_{ijkl}$, and $\breve{\eta}_{ijkl}(\mathbf{x})=(c \,/\, \beta_l)\, \eta_{ijkl}$, while the rest of the coefficients are as defined in elastic PML formulation as given in equation \eqref{eq:pmlElastic}.
\end{appendices}

\end{document}